\documentclass[10pt]{amsart}
	\usepackage[all]{xy}
	\usepackage{amssymb}
	\usepackage{graphics}

	\newcommand{\Spe}{{\rm Spec}}

	\newcommand{\T}{{\mathcal{T}}}

	\renewcommand{\O}{{\mathcal{O}}}
	\newcommand{\Z}{{\mathbb{Z}}}
	\newcommand{\lc}{\left\lceil}
	\newcommand{\rc}{\right\rceil}

	\newcommand{\be}{\begin{equation}}
	\newcommand{\ee}{\end{equation}}

	\newtheorem{pro}{Proposition}[section]
	\newtheorem{lemma}[pro]{Lemma}
	
	\newtheorem{cor}[pro]{Corollary}
	\newtheorem{theorem}[pro]{Theorem}

\begin{document}
	\bibliographystyle{amsplain}
	\title{A homological approach for computing the 
	tangent space of the deformation functor of curves with 
	automorphisms. 
	}

	\author{A. Kontogeorgis}

	\begin{abstract}
		We give an alternative approach to the computation 
		of the dimension of the tangent space of the 
		deformation space of curves with automorphisms. 
		A homological version of the local-global principle similar 
		to the one of J.Bertin, A. M\'ezard is proved, and a computation 
		in the case of ordinary curves is obtained, by application 
		of the results of S. Nakajima for the Galois module structure
		of the space of 2-holomorphic differentials on them. 
        \end{abstract}

	\email{kontogar@aegean.gr}
	\address{
	Department of Mathematics, University of the \AE gean, 83200 Karlovasi, Samos,
	Greece\\ { \texttt{\upshape http://eloris.samos.aegean.gr}}
	}
	\date{\today}

	\maketitle

	\section{Introduction}

	Let $X$ be a  non-singular curve of genus $g\geq 2$ 
	defined over an algebraic closed field 
	of positive characteristic, together with a subgroup $G$ of the  automoprhism
	group. 
	In \cite{Be-Me} J.Bertin, A. M\'ezard 
	proved that the equivariant cohomology of Grothendieck 
	$H^1(G,\T_X)$ measures the tangent space of the global deformation functor
	of curves with automorphisms. 
	This  dimension is a measure in how many directions a curve can be deformed
	together with  a subgroup of an automoprhism group. 

	Since  the genus $g$ of $X$ is $g\geq 2$  the edge homomorphisms of the 
	spectral sequence of Grothendieck \cite[5.2.7]{GroTo} give as that 
	\[
	H^1(G,\T_X) =H^1(X,\T_X)^G, 
	\]
	{\em i.e.} the first  equivariant cohomology equals the 
	$G$-invariant space of the ``tangent space" of the moduli space at $X$. 
	J. Bertin and A. M\'ezard were successful to compute $H^1(G,\T_X)$, 
	using a version of equivariant Chech theory, and prove a local-global
	theorem. 

	It is tempting, in order to compute the $G$-invariants of the 
	space $H^1(G,\T_X)$ to use Serre's duality to pass to the 
	space $H^0(G,\Omega^{\otimes 2}_X)$. One should be careful using 
	this approach,
	because since we are  considering  the dual space of $H^1(X,\T_X)$ 
	it is not the functor of invariants that we have to 
	consider, but the 
	  the adjoint functor, {\em i.e.}, the functor 
	of covariants.

	This article consists of two parts. In first part, we 
	use the normal basis 
	theorem for Galois extensions and the explicit form 
	of Serre duality in terms of repartitions 
	\cite[7.14.2]{Hartshorne:77},\cite[I.5]{StiBo}, 
	in order to compute the covariant elements of
	$H^0(X,\Omega_X^{\otimes 2})$. This leads us to  
	a homological version 
	of the local-global theorem of J. Bertin A. M\'ezard (\ref{mainprop}).
	The duality between the homological and 
	the cohomological approach is 
        emphasized by defining  a cap product 
	\[
	H^p (G, H^1(G,\T_X) ) \times H_q (G, H^0(X,\Omega_X^{\otimes 2} )) \rightarrow 
	H_{q-p} (G,k).
	\]

	In second part we try to apply known results on 
	the Galois module structure of the space of 
	holomorphic differentials, in order to compute covariant 
	elements. 
	We have to notice here, that, as far the author knows, 
	if the characteristic $p$ divides $|G|$, the Galois module 
	structure of $H^0(X,\Omega^{\otimes s})$ is far from 
	being understood, and there are only partial results 
	mainly in the case of tame ramification 
	\cite{Kani:88},\cite{Nak:87i},\cite{Nak:85j}, or in 
	the case of ordinary curves. 

	More precisely, 
        S. Nakajima in \cite{Nak:85} studied the Galois module 
	structure of the ``semi-simple part" of $ H^0(X,\Omega(-D)$
	with respect to the Cartier operator if $G$ is a $p$-group and 
	$D$ is an effective $G$-invariant divisor on $X$. 
        Thus, for the Zariski dense set in the moduli space of curves 
	of genus $g$,  of ordinary curves (curves 
		the Galois module structure of $H^0(X,\Omega^{\otimes 2})$
	is known.
        
	The tangent space to the deformation functor of ordinary curves 
	was studied by G. Cornelissen and F. Kato in \cite{CK}.
	We obtain a weaker result than their 
	result by using the methods of S. Nakajima.

		\section{Computations}

	Denote by $K_X$ the function field of the curve $X$. 
	Let $\mathcal{K}_X$ be   the constant sheaf $K_X$, 
	and consider the exact sequence of sheaves
	\begin{equation} \label{repart1}
	0 \rightarrow \O_X \rightarrow \mathcal{K}_X \rightarrow 
	\frac{\mathcal{K}_X}{\O_X} \rightarrow 0. 
	\end{equation}
	The sheaf $\frac{\mathcal{K}_X}{\O_X}$ can be expressed in the 
	form 
	\[
	\frac{\mathcal{K}_X}{\O_X}=\bigoplus_{P\in X} i_*(K_X/\O_P),
	\]
	where $i:\Spe \O_P \rightarrow X$ is the inclusion map. 

	We tensor the  sequence (\ref{repart1}) with the sheaf $\Omega_X^{\otimes 2}$ over 
	$\O_X$
	and get the sequence:
	\[
	0 \rightarrow \Omega_X^{\otimes 2}  \rightarrow \mathcal{K}_X 
	\otimes _{\O_X} \Omega_X^{\otimes 2}  \rightarrow 
	\bigoplus_{P\in X} i_*(K_X/\O_P)\otimes \Omega_X^{\otimes 2}  \rightarrow 0. 
	\]
	We will denote by  $\mathcal{M}^{\otimes 2}= \mathcal{K}_X 
	\otimes _{\O_X} \Omega_X^{\otimes 2} $ the sheaf of meromorphic 2-differentials
	and by $\Omega_P^{\otimes 2}=\Omega_X^{\otimes 2} \otimes _{\O_X} \O_P$.
	Thus we might write 
	\[
	\bigoplus_{P\in X} i_*(K_X/\O_P)\otimes \Omega_X^{\otimes 2} =
	\bigoplus_{P\in X} i_*(\mathcal{M}^{\otimes 2} /\Omega_P^{\otimes 2}).
	\]

	We apply the global section functor:
	\[
	0\rightarrow  \Gamma(X,\Omega_X^{\otimes 2})  \rightarrow \Gamma(X,\mathcal{K}_X 
	\otimes _{\O_X} \Omega_X^{\otimes 2})  \rightarrow 
	\bigoplus_{P\in X} i_*(\mathcal{M}^{\otimes 2} /\Omega_P^{\otimes 2})\rightarrow 
	H^1(X, \Omega_X^{\otimes 2}) \rightarrow \cdots
	\]
	Since $X$ is a curve of genus $g\geq 2$ we have that
	$H^1(X, \Omega_X^{\otimes 2})=0$ and if we denote by $\Omega=\Gamma(X,\Omega_X^{\otimes 2}) $
	and $M=\Gamma(X,\mathcal{M}^{\otimes 2})$ 
	the spaces of global sections 
	of homomorphic and meromorphic differentials we have:
	\begin{equation} \label{exact-dif}
	0 \rightarrow \Omega \rightarrow M \rightarrow
	\Gamma \left(X,\bigoplus_{P\in X} 
	i_*(\mathcal{M}^{\otimes 2}/\Omega_P^{\otimes 2})\right)
	 \rightarrow 0.
	\end{equation}

	\begin{lemma}
		The $G$-module $M$ as a  $K_Y[G]$ module is projective.
	\end{lemma}
	\begin{proof}
	Let $w$ be a meromorphic differential of the curve $Y=X/G$,
	and denote by $K_Y$ the function field of the curve $Y$. The lift 
	$\pi_* w$ is a $G$-invariant meromorphic differential on $X$, 
	and $M$ can be recovered as the set of the expressions
	\[
	M=\{f \cdot \pi_*(w), \;\;\; f \in K_X\}.
	\]
	We want to apply the functor of covariants, {\em i.e.}, 
	to tensor with $K_Y \otimes_{K_Y[G]}$. 
	We notice first that by the normal basis theorem \cite[6.3.7 p.173]{Weibel}
	for the Galois
	extension $K_X/K_Y$ we obtain that $K_X \cong  K_Y[G]$ as a Galois module, thus 
	$M$ is isomorphic to $K_Y[G]$ as a $K_Y[G]$-module and the desired
	result follows. 
	\end{proof}
	We consider the long exact homology sequence arising from 
	(\ref{exact-dif}) after taking the functor of covariants:
	\begin{equation} \ref{mainprop}
	\cdots \rightarrow H_1(G,M) \rightarrow 
	H_1(G,\bigoplus_{P\in X} i_*(\mathcal{M}^{\otimes 2}/\Omega^{\otimes 2}_P)
	 )
	\rightarrow 
	\Omega_G \stackrel{\alpha}{\rightarrow}
	M_G \rightarrow \end{equation}
	\[
	\rightarrow
	\Gamma(X,\bigoplus_{P\in X} i_*(\mathcal{M}^{\otimes 2}/\Omega^{\otimes 2}_P)
	)_G 
	\rightarrow 0.
	\]
	Since $M \cong K_Y[G]$ we have  $ H_1(G,M)=0$ and $M_G=\{ f\cdot \pi_*(w)\}$, 
	with $f\in K_Y$.

	Thus 
\[
\Omega_G= H_1(G,
\Gamma(X,\bigoplus_{P\in X} i_*(\mathcal{M}^{\otimes 2}/\Omega_P^{\otimes 2})
 ))
\oplus \mathrm{Im} \alpha.
\]
{\bf Remark:} If the order $|G|$ of the group $G$ is  prime to  the
characteristic  $p$ then 
the order $|G|$
is
invertible in the module 
$\Gamma(X,\bigoplus_{P\in X} i_*(\mathcal{M}^{\otimes 2}/\Omega_P
^{\otimes 2}))
 $ and 
the first homology is zero, therefore 
\[
\Omega_G = \mathrm{Im} \alpha.
\]
\begin{pro}
	Let $b_1,\ldots,b_r$ be the set of ramification points of the 
	cover $X\rightarrow Y$, and let $G_i=G(b_i)$ be the 
	corresponding decomposition groups. The following holds:
	\[
	H_1(G,\bigoplus_{P\in X} i_*(\mathcal{M}^{\otimes 2} /
	\Omega_P^{\otimes 2}))=
	\bigoplus_{i=1}^r 
	H_1(G_i, \mathcal{M}^{\otimes 2}/\Omega^{\otimes 2}_{b_i})).
	\]
\end{pro}
\begin{proof}
	Let $P$ be a point  of $X$, and let $t_P$ be a local uniformizer
	at the point $P$. Consider an element 
	$a=\sum_{P\in X} a_P P \in 
	\bigoplus_{P\in X} i_*(\mathcal{M}^{\otimes 2}/\Omega^{\otimes 2}_P))
	$. 
	The polar part of $a$ at $P$ is equal to 
	$\sum_{\nu=-n}^{-1} \frac{a_\nu}{t^\nu}$. For an element 
	$g\in G$ we have that 
	\[
	g\left(\sum_{\nu=-n}^{-1} \frac{a_\nu}{t^\nu}\right)=
	\sum_{\nu=-n}^{-1} \frac{a_\nu}{g(t)^\nu}
	.\]
	The element $g(t)$ is the local uniformizer at the point 
	$g(P)$. 
This proves that the action of the element $g\in G$ on $a$ is of the form 
\[
g(\sum_{P\in X} a_P P)= \sum _{P \in X} a_P g(P).
\]
Let $M_P= i_*(\mathcal{M}^{\otimes 2}/\Omega_P^{\otimes2}) $
be the summand corresponding to the point $P$, and let  
 $G(P)=\{g\in G: g(P)=P\}$ be the decomposition group at the point $P$. 
 We consider the induced module, seen as a subspace of $\oplus_{P\in X} 
 \mathcal{M}^{\otimes 2}/\Omega^{\otimes 2}_P$,
 \[
 \mathrm{Ind}_{G(P)}^G M_P = K_Y[G] \otimes _{K_Y[G(P)]} M_P =
 \bigoplus_{g\in G/G(P)} M_{g(P)} .   
 \]
 Shapiro's lemma \cite[6.3.2]{Weibel} implies that 
 \[
 H_1(G,\mathrm{Ind}_{G(P)}^G M_P)=H_1(G(P),M_P).
 \]
Thus if $P$ is not a ramification point it does not contribute to 
the cohomology, and the desired formula comes by the sum of the contributions 
of the  ramification groups. 

\end{proof}
We have proved that the following sequence is exact:
\begin{equation} \label{mainprop}
	0 \rightarrow \bigoplus_{i=1}^r 
	H_1(G(b_i), \mathcal{M}^{\otimes 2}/\Omega^{\otimes 2}_{b_i})) \rightarrow 
	\Omega_G \rightarrow \mathrm{Im}\alpha \rightarrow 0 
\end{equation}
which is exactly the dual sequence of J. Bertin, A. M\'ezard 
\cite[p . 206]{Be-Me}.
\begin{pro}
	Let $b_1,\ldots, b_r$ be the ramification points of the cover 
$\pi:X \rightarrow Y$, and assume that the groups in the ramification 
filtration at each ramification point $b_k$ have orders
\[
e_0^{(k)} \geq e_1^{(k)}\geq \cdots \geq e_{n_k}^{(k)}>1. 
\]
The dimension of the space $\mathrm{Im} \alpha$ is given by:
\[
\dim_k \mathrm{Im} \alpha= 3g_y-3 + \sum_{k=1}^r \lc  2\sum_{i=1}^{n_k} 
\frac{e_i^{(k)}-1}{e_0^{(k)} }
\rc .
\]
\end{pro}
\begin{proof}
We are looking for elements of the form $f \pi^*(w)$, $f\in K_Y$, 
such that $\mathrm{div}_X f \pi^*(w) \geq 0$. 
We know that if $w$ is a $2$-differential then:
\[
\mathrm{div} (\pi^* (w)) =\pi^* (\mathrm{w}) +2 R,
\]
where $R$ is the ramification divisor of $\pi:X \rightarrow Y$. 
Therefore $\pi^*(w)$ is holomorphic if and only if 
$\pi^* (\mathrm{w}+2R) \geq 0$.

We will push forward again and we will use Riemann-Roch on $Y$. 
We want to compute the dimension of the space 
\[
L(K+ \pi_*(2R)/|G|),
\]
where $K$ is the canonical divisor. 

The ramification divisor is 
\[
R=\sum_{k=1}^r \sum_{i=1}^{n_k} ( e_i^{(k)}-1) b_k,
\]
and  Riemann-Roch theorem implies that 
\[
\dim_k \mathrm{Im} \alpha= 3g_y-3 + \sum_{k=1}^r \lc  2\sum_{i=1}^{n_k} 
\frac{e_i^{(k)}-1}{e_0^{(k)} }
\rc .
\]
\end{proof}
\begin{pro}
	For a ramification point $b$ we have $H_1 (G(b),M_b ) \cong V(b)_G$, 
	where $V(b)$ is a finite $k$-vector space, with known $G$-module structure. 
\end{pro}
\begin{proof}
	The space $M_b$ consists of elements in $1/t k[1/t]$.
	Let $f(1/t)\in k[1/t]$ be a $G$-invariant element. Then we can 
	consider the direct product of $G$-modules:
	\[
	k((t)) = 1/tf(1/t) k(1/t) \bigoplus L( v(f)),
	\]
	where $k((t))$ is the field of formal Laurent series, 
	$v$ is the  valuation of the field $k((t))$, 
	and $L\big(v(f)\big)=\{g \in k((t)), v(g)+v(f(1/t)) \geq 0\}$.
	Theorem 90 of Hilbert implies that $H_1\big(G(b),k((t)) \big)=0$, 
	and since 
	\[H_1\big(G(b),k((t)) \big) =H_1 \big(G(b), f(1/t) k(1/t)) \big) \bigoplus
	H_1\big(G(b),L( v(f)) \big),\]
	we have that $H_1\big(G(b), 1/tf(1/t) k(1/t)) \big) =0$. 
	We now consider the short exact sequence:
	\[
	0 \rightarrow 1/t  k[1/t] \stackrel{ \times f(1/t) }{
         -\!\!\!-  \!\!\! - \!\!\! \longrightarrow} f(1/t) 1/t  k[1/t] 
	\rightarrow V(b) \rightarrow 0,
	\]
	and by a dimension shifting argument we obtain the desired result. 
\end{proof}

\section{Definition of a cap product.}
In this section we will define a cap product between the homology 
and cohomology groups of $A:=H^1(X,\T_X)$ and $B:=H^0(X,\Omega_X^{\otimes 2})$. 
Serre duality implies the existence of a trace function
\[
A \times B  \rightarrow k.
\]
Thus, a cap product can be defined (\cite[p.113]{BrownCoh}) 
\[
H^p(G, A)  \otimes H_q(G, B) \rightarrow
H_{q-p} (G, A\otimes B) 
\stackrel{ \mathrm{tr} }{\longrightarrow} H_{q-p}(G,k).
\]
In particular,
the above cap product gives us the pairing
\[
H^1(G,A) \otimes H_1( G,B) \rightarrow k
\]
connecting the homological and cohomological approaches to 
the theory. 

{\bf Remark:}
Since $G$ is a finite group, the Tate cohomology groups are defined, 
and they seem a more  natural tool for the study  
	of cap products and duality. Unfortunately we are 
	interested for the computation of invariants and co-invariants, 
	{\em i.e.} for low index cohomology groups, and Tate cohomology 
	can not be applied here. 
\section{Ordinary Curves}
Let $D$ be an effective divisor on the curve $X$, we will denote by 
\[
\Omega_X(-D)=\{ f \in k(X): \mathrm{div}(f) + \mathrm{div}(\omega) \geq -D\}.
\]
S. Nakajima in \cite{Nak:85} provided us with a method for  
computing the Galois module structure of the semisimple part of 
$\Omega_X(-D)$ with respect to the Cartier operator on $X$. 

If the curve $X$ is ordinary then the semisimple part of the 
Cartier operator is identified with the space $\Omega_X(-D)$ itself.
We are interested in computing the space of covariants
$\Omega^{\otimes 2}_G$ of the holomorphic $2$-differentials.

The space of holomorphic $2$-differentials can be identified 
with the space 
\[
\{f \in k(X) : 
\mathrm{div}(f) + 2\mathrm{div}(\omega) \geq 0\} =\{ \mathrm{div}(f)+ 
\mathrm{div}(\omega) \geq - \mathrm{div}(\omega) \}
\]

\begin{lemma}
	There is a $G$-invariant  differential $\omega$ in $X$, 
	such that $\mathrm{div}(\omega)$ is effective, and has support 
	that does not intersect the branch locus. 
\end{lemma}
\begin{proof}
	Let $b_1,\ldots,b_r$ be the ramification points of the cover 
	$\pi:X \rightarrow Y=X/G$. 

	Let $\phi_1$ be an arbitrary meromorphic differential on $Y$.
	We will select a meromorphic differential $\phi=f \phi_1$ on the curve 
	$Y$ such that 
	$\mathrm{div}(\phi)=\mathrm{div}(f\phi_1) + A \geq 0$,
	where $A$ is the divisor
	\[
	A:= \sum_{i=1}^r \lc \sum_{i=0}^{\infty} \frac{e_i(b_i) -1}{e_0(b_i)}
	\rc b_i.
	\]

	Notice that if we assume that we are working on an ordinary curve 
	then $e_2(b_i)=0$ \cite{Nak}, and the above divisor can be written 
	as 
	\[
	A=\sum_{i=1}^r \lambda_i b_i,
	\]
	where $\lambda_i=1$ if $b_i$ is ramified tamely and $\lambda_i=2$ if
	$b_i$ is ramified wildly.
	
	This means that we are looking for a function 
	$f\in L_Y(K +A)$. 
	Using Riemann-Roch we compute 
	\[
	\ell(K+A)= \ell(-A) + 2g_Y-2 + \deg(A)-g_Y +1=g_Y -1 +r+s,
	\]
	where $s$ is the number of wild ramified branch points. 
	For such a selection of $f$ we have that 
	\[
	\mathrm{div}( \pi^* (f \phi_1)=\pi^* \mathrm{div}(f\phi_1) +R \geq 0,
	\]
	where $R$ is the ramification divisor given by 
	\[
	R=\sum_{i=1}^r  \sum_{P \mapsto b_i} 
	\sum_{i=0}^{\infty} \big({e_i(P) -1}\big).
	\]
	Moreover the divisor $\pi^* (f \phi_1)$ is $G$ invariant, and we 
	can select $f \in L_Y(K +A)$ such that is has polar divisor $A$. 
	This imply that the support of $\pi^* (f \phi_1)$ has no intersection 
	with the branch locus. 
\end{proof}
	
We can now apply the method of S. Nakajima on 
$\Omega_X^{\otimes 2}=\Omega_X( \pi^*(f \phi_1))$. Let $S$ be the set of points
of the curve $Y$ such that 
$\pi^{-1}(S)=\mathrm{supp} \big(\mathrm{div} (\pi^*(f \phi_1))\big)$. We follow 
the notation of \cite{Nak:85}. Let 
$S_0=\{b_1,\ldots,b_r\}$. For each $i=1,\ldots,r$ we choose a point $P_i\in X$, 
satisfying $\pi(P_i)=b_i$, and let $G_i$ be the decomposition group 
at $P_i$. We consider the $k[G]$-modules 
$k[G/G_i]=\{\sum_{\sigma \in G/G_i} a_\sigma \sigma \}$. We define surjective 
$k[G]$-homomorphisms $\Phi_i:k[G/G_i]\rightarrow k$, by 
$\Phi_i(\sum_{\sigma\in G/G_i} a_\sigma \sigma)=\sum_{\sigma \in G/G_i} a_\sigma$,
and also 
$\Phi( (\xi_1,\ldots,\xi_r))=\sum_{i=1}^r \Phi_i (\xi_i)$.
\begin{theorem}
	Let $G$ be a $p$-group.
	The Galois module structure of $\Omega_X^{\otimes 2}$ is determined 
	by the following exact sequence:
	\begin{equation}\label{gal-struct}
	0 \rightarrow \Omega_X\big( - \pi^{-1}(S) \big) \rightarrow 
	\Omega_X(-\pi^{-1}\big(S \cup S_0)\big) \rightarrow \ker\Phi \rightarrow 0,
\end{equation}
	where $\Omega_X(-\pi^{-1}\big(S \cup S_0)\big)\cong k[G]^{3g_Y-3 +2r}$.
\end{theorem}
\begin{proof}
	Following the method of Nakajima, 
we define a $k[G]$-homomorphism 
\[
\phi: \Omega_X(-\pi^{-1}\big(S \cup S_0)\big) \rightarrow 
\bigoplus_{i=1}^r k[G/G_i], 	
\]	
by
\[
\phi (\omega) =\bigoplus_{i=1}^r \left(
\sum_{\sigma \in G/G_i} \mathrm{Res}_{\sigma P_i} \sigma
\right), \quad \omega\in \Omega_X(-\pi^{-1}\big(S \cup S_0)\big).
\]
Then $\ker(\phi)=\Omega_X\big( - \pi^{-1}(S) \big)$,
and $\mathrm{Im}(\phi)=\ker \Phi$, using the residue and Riemann-Roch 
theorems. 

	For the $k[G]$-structure of $\Omega_X(-\pi^{-1}\big(S \cup S_0)\big)$
 Theorem 1 of \cite{Nak:85} implies that it is a free $k[G]$-module of rank 
 $\gamma_Y-1 + |S|+|S_0|$, where $\gamma_Y$ is the $p$-rank of the Jacobian 
 of $Y$. Since the curve is ordinary we have $g_Y=\gamma_Y$. On the 
 other hand $|S|=2g_Y-2+r$ and $|S_0|=r$. 
\end{proof}

The module in the middle of equation (\ref{gal-struct}) is $k[G]$-projective, therefore
it implies the following long exact sequence:
\[
0 \rightarrow H_1 (G, \ker \Phi) \rightarrow 
\left(\Omega_X^{ \otimes 2} \right) _G \rightarrow
 k[G]^{3g_Y-3 +2r}_G \rightarrow \ker \Phi_G \rightarrow 0.
\]
This implies that the desired dimension can be computed: 
\begin{equation} \label{comp:dim1}
\dim_k  \left(\Omega_X^{ \otimes 2} \right) _G 
=\dim_k H_1(G, \ker \Phi) + 3g_Y-3 +2r -\dim_k \ker \Phi_G.
\end{equation}
We will use the sequence 
\begin{equation} \label{sec:ses}
0\rightarrow \ker \Phi \rightarrow \bigoplus_{i=1}^r k[G/G_i] 
\stackrel{\Phi}{\longrightarrow} k \rightarrow 0,
\end{equation}
in order to compute the homology groups of $\ker \Phi$.

Equation  (\ref{sec:ses}), gives  the
long exact sequence:
\begin{equation} \label{long:res1}
H_2(G , \bigoplus_{i=1}^r k[G/G_i])\stackrel{\psi_1}{\longrightarrow}  
H_2(G,k) 
 \rightarrow 
H_1(G,\ker \Phi) \rightarrow H_1(G,\bigoplus_{i=1}^r k[G/G_i]) 
\stackrel{\psi_2}{\longrightarrow}
\end{equation}
\[
\stackrel{\psi_2}{\longrightarrow}
  H_1(G,k) \rightarrow \ker \Phi _G \rightarrow 
\left(\bigoplus_{i=1}^r k[G/G_i] \right)_G \rightarrow k_G \rightarrow 1.
\]
Using the above sequence we compute:
\[
\dim_k H^1( G, \ker \Phi) = \dim_k \mathrm{Coker}\psi_1 + \dim_k \ker \psi_2.
\]
It is known that $H_1(G,k)=\frac{G}{[G,G]}\otimes _\Z k$
and Hopf's theorem \cite[6.8.8]{Weibel} implies 
\[
 H_2(G,k)= \frac{R \cap [F,F]}{[F,R]}\otimes_\Z k,
\]
where $1 \rightarrow R \rightarrow F \rightarrow G \rightarrow 1$ 
is a free presentation of $G$.

Shapiro's lemma \cite[p. 73]{BrownCoh},  gives  that 
$H_1(G,k[G/G_i]) = H_1(G_i,k)$. Moreover, since the curve in 
question is ordinary the ramification groups are elementary 
abelian, thus $H_1(G,k[G/G_i])=G_i\otimes_\Z k$. 
Therefore, for the computation of the kernel of $\psi_2$ we have:
\[
\psi_2: \bigoplus_{i=1}^r G_i\otimes _\Z k \rightarrow 
\frac{G}{[G,G]}\otimes _\Z k.
\]
The kernel of $\psi_2$ equals $\cap_{i=1}^r (G_i \cap [G,G])\otimes_\Z k$. 
This is a group theoretic description of the kernel of $\psi_2$.

Using equation
(\ref{long:res1}) we compute:
\begin{equation} \label{lala2}
\dim_k H_1(G,\ker \Phi)= \dim_k \mathrm{coker}(\psi_1) +
\sum_{i=r}^r  \log_p|G_i| - 
\dim_k \frac{G}{[G,G]}\otimes _Z k + \dim_k \ker\Phi_G -r +1.
\end{equation}
Combining (\ref{comp:dim1}) with (\ref{lala2}) we obtain
\[
\dim_k (\Omega_X^{\otimes 2})_G =\dim_k \mathrm{coker}(\psi_1) 
+ 3g_Y-3 +r +
\sum_{i=1}^r  \log_p |G_i| - 
\dim_k  \frac{G}{[G,G]}\otimes _\Z k   +1.
\]

For the computation of the cokernel of $\psi_1$ we proceed as follows:
For every ramification group fix a set of generators $F_i \in F$ and
consider a set of relations $R_i$, such that $G_i=F_i/R_i$. 
Using   Hopf's theorem for the computation of $H_2(G,k)$, the 
study of the map $\psi_1$ is reduced to the study of the map:
\[
\bigoplus_{i=1}^r  \frac{R_i \cap [F_i,F_i]}{[F_i,R_i]} \rightarrow 
\frac{R \cap[F,F]}{[F,R]}. 
\]

The groups $G_i$ are elementary abelian. If $G_i$ is a cyclic group of 
order $p$ then it is immediate from Hopf's theorem that 
$H_2(G_i,k)=0$. 

We consider now the case of groups $G_i$ that have at least 
two cyclic summands. 
Since the groups $G_i$ are  abelian, we have $[F_i,F_i] \subset
R_i$, thus $R_i \cap [F_i,F_i]=[F_i,F_i]$. 
For a given $i=1,\ldots,r$ consider the map 
\[f_i:
\frac{[F_i,F_i]} {[F_i,R_i]} \rightarrow \frac{R \cap [F,F]}{[F,R]}.
\]
The kernel of $f_i$ is $(\ker f_i=[F_i,F_i] \cap [F,R])/[F_i,R_i]$ and the 
image is isomorphic to 
\[
\mathrm{Im}(f_i) \cong \frac{[F_i,F_i]}{[F_i,F_i] \cap [F,R]}
\otimes _\Z k
\cong
\frac{[F_i,F_i][F,R]}{[F,R]}
\]
Combining this information together for all $i$ such that 
$\log_p |G_i| >1$ we obtain:
\[
\mathrm{coker}{\psi_1}=
\frac{R \cap [F,F]}{\langle  [F_i,F_i][F,R]
\rangle} \otimes_\Z k.
\]
where $i$ runs over the  ramification points such that $\log_p |G_i| \geq 1$. 
If all $G_i$ have order $p$ then the above formula reduces to Hopf's formula 
for $H_2(G)$. We collect all pieces of computation in the following 
\begin{pro}
Let $G$ be a $p$-group that is a subgroup of the automorphism 
	group of an ordinary curve, and let $G_i$ be the decomposition 
	groups at the ramification points. 	Using the above notation 
	we have for the dimension of covariant 2-differentials
	
	\[
	\dim_k \left(\Omega_X^{\otimes 2}\right)_G =
	3g_Y -3 +r + \sum_{i=1}^r \log_p |G_i| -
	\dim_k \frac{G}{[G,G]}\otimes_\Z k -1 +
	\dim_k \frac{R \cap [F,F]}{\langle  [F_i,F_i][F,R]
	\rangle}  \otimes_\Z k.
	\]
\end{pro}

Comparison of the computation done so far with the result of 
G.Cornelissen, F.Kato implies the following corollary:
\begin{cor}
	Let $G$ be a $p$-group that is a subgroup of the automorphism 
	group of an ordinary curve, and let $G_i$ be the decomposition 
	groups at the ramification points. Suppose that $p>3$. 
	Using the above notation we 
	obtain:
	\[
	\dim_k \frac{R \cap [F,F]}{\langle  [F_i,F_i][F,R]
	\rangle} \otimes_\Z k = \dim_k \frac{G}{[G,G]} \otimes_\Z k -1.
	\]
\end{cor}

\providecommand{\bysame}{\leavevmode\hbox to3em{\hrulefill}\thinspace}
\providecommand{\MR}{\relax\ifhmode\unskip\space\fi MR }
\providecommand{\MRhref}[2]{%
  \href{http://www.ams.org/mathscinet-getitem?mr=#1}{#2}
}
\providecommand{\href}[2]{#2}

\end{document}